\documentclass[11pt]{amsart}
\usepackage{mathrsfs}
\usepackage{amsfonts}
\usepackage{epsfig}
\usepackage{amsmath,amssymb,amsthm}

\newcommand{\co}{\mathbb{C}}

\newcommand{\rp}{\mbox{\rm Re}}

\newcommand{\p}{\partial}
\newcommand{\g}{\lambda}
\newcommand{\ol}{\overline}

\begin{document}

\newtheorem{thm}{Theorem}[section]
\newtheorem{prop}[thm]{Proposition}
\newtheorem{lem}[thm]{Lemma}
\newtheorem{cor}[thm]{Corollary}
\newtheorem{conj}[thm]{Open Problem}

\theoremstyle{definition}

\newtheorem{defn}[thm]{Definition}
\newtheorem{notation}[thm]{Notation}
\newtheorem{example}[thm]{Example}

\numberwithin{equation}{section}

\title[ Positive invariant harmonic function in the unit ball]
{On the monotonicity of positive invariant harmonic functions in the
  unit ball} 
\author{Yifei Pan}
\address{Yifei Pan\\
Department of Mathematical Sciences \\ Indiana University -
Purdue University Fort Wayne \\ Fort Wayne, IN 46805-1499}
\address
{School of Mathematics and Informatics \\ 
Jiangxi Normal University, Nanchang, China}
\email{pan@ipfw.edu}

\author{Mei Wang}
\address{Mei Wang\\
Department of Statistics\\ University of Chicago\\
Chicago, IL 60637}
\email{meiwang@galton.uchicago.edu}


\begin{abstract}
  A monotonicity property of Harnack inequality is proved for
  positive invariant harmonic functions in the unit ball.
\end{abstract}

\maketitle

%
\section{Introduction}\label{sec0}
%

Let $B^n=\{x\in\mathbb{R}^{n}:|x|<1\},~n\ge 2$ be the open unit ball in
$\mathbb{R}^n$. $S^{n-1}=\partial B^n$.  Consider the differential operator
$$\Delta_{\lambda}=(1-|x|^2)\left\{\frac{1-|x|^2}{4}\sum_{j} \frac{\p^2}{\p
    x_j^2} +\g\sum_{j}x_j\frac{\p}{\p
    x_j}+\g(\frac{n}{2}-1-\g)\right\},~\g\in\mathbb{R}.$$
In this paper, we prove a monotonicity property of invariant harmonic
functions that are solutions of $\Delta_{\lambda}u=0$ and are defined by
positive Borel measures on the sphere with respect to the Poisson
kernel $P_\lambda$ (see below).  

\vspace{3mm} This section describes the theorems and their
corollaries.  The proofs are provided in the next two sections.  

          %
          %
%
\begin{thm}  
  Let $u$ be a positive invariant harmonic function defined in $B^n$
  by a positive Borel measure $\mu$ on $S^{n-1}$ with the Poisson
  kernel $P_\lambda$.  For $\zeta\in S^{n-1}$, if $\g>-\frac{n}{2}$
  (if $\g<-\frac{n}{2})$, the function
$$\dfrac{(1-r)^{n-1}}{(1+r)^{1+2\g}}u(r\zeta)$$
is decreasing (increasing) for $0\leq r<1$, and the function
$$\dfrac{(1+r)^{n-1}}{(1-r)^{1+2\g}}u(r\zeta)$$
is increasing (decreasing) for $0\leq r<1$.  Also
$$\lim_{r\to 1}(1-r)^{n-1}u(r\zeta)=
\begin{cases}
2^{1+2\g}\mu(\{\zeta\}), 
& \lambda > -\frac{n}{2}\\
\infty, & \lambda<-\frac{n}{2}, ~ \mu(\{\zeta\}^c)>0\\
2^{1+2\g}\mu(\{\zeta\}), & \lambda<-\frac{n}{2}, ~ \mu(\{\zeta\}^c)=0
\end{cases}
$$
and 
$$\lim_{r\to
  1}\dfrac{u(r\zeta)}{(1-r)^{1+2\g}} =
\int_{S^{n-1}}\dfrac{2^{1+2\g}}{|\zeta-\xi|^{n+2\g}}d\mu(\xi).$$
\end{thm}
%

          %
          %
%
\vspace{3mm}\noindent{\bf Remarks.} ~ 

\begin{enumerate}
\item Invariant harmonic functions are the solutions of $\Delta_\lambda u=0$.
These solutions also satisfy certain invariance property with respect
to M\"obius transformation.  Invariant harmonic functions generally do
not possess good boundary regularity, as shown in Liu and Peng
{\cite{LiuPeng}}. 

\vspace{1mm} \item Let $\mu$ be a positive Borel measure on $S^{n-1}$ and
$P_\g$ be the Poisson kernel
$$P_\g(x, \zeta)=\frac{(1-|x|^2)^{1+2\g}}{|x-\zeta|^{n+2\g}}.$$
It is known that the integral 
$$u(x)=\int_{S^{n-1}}P_\g(x,\zeta)d\mu(\zeta)$$
defines an invariant harmonic function in $B^n$ (\cite{axler}, p. 119).

\vspace{1mm}\item On the completion of the current work, we learned
that the limit cases for $n=2, \lambda=0$ in Theorem 1.1 were obtained
by Simon and Wolff (\cite{SimonWolff}, ref. Chapter 10,
p. 546 in \cite{simon}).

\vspace{1mm}\item The critical value $\lambda=-\frac{n}{2}$ yields
the degenerate case with the constant Poisson kernel.

\end{enumerate}

\vspace{3mm}The following theorem characterizes the behavior of
invariant harmonic functions on the rays.

          %
          %
%
\begin{thm}  
  Let $u$ be a positive invariant harmonic function defined in $B^n$
  by a positive Borel measure $\mu$ on $S^{n-1}$ with the Poisson
  kernel $P_\lambda$.  Let $\zeta\in S^{n-1}$ and $~0\le r'\le r<1$.

\vspace{2mm}\noindent If $\g>-\frac{n}{2}$,
\begin{eqnarray*}
\left(\!\frac{1\!-r}
{1\!-r'}\!\right)^{\!\!2\lambda+1}
\left(\!\!\frac{1\!+\!r'}{1\!+\!r}\right)^{\!\!n-1} 
\!\!\!\!\!u(r'\zeta)
\le u(r\zeta) \le
\left(\!\frac{1\!+r}
{1\!+r'}\!\right)^{\!\!2\lambda+1}
\left(\!\!\frac{1\!-\!r'}{1\!-\!r}\right)^{\!\!n-1} 
\!\!\!u(r'\zeta)
\end{eqnarray*}

\vspace{2mm}\noindent If $\g<-\frac{n}{2}$, 
\begin{eqnarray*}
\left(\!\frac{1\!+r}
{1\!+r'}\!\right)^{\!\!2\lambda+1}
\left(\!\!\frac{1\!-\!r'}{1\!-\!r}\right)^{\!\!n-1} 
\!\!\!\!\!u(r'\zeta)
\le u(r\zeta) \le
\left(\!\frac{1\!-r}
{1\!-r'}\!\right)^{\!\!2\lambda+1}
\left(\!\!\frac{1\!+\!r'}{1\!+\!r}\right)^{\!\!n-1} 
\!\!\!u(r'\zeta)
\end{eqnarray*}

\vspace{2mm}\noindent For $r'=0$, the above becomes
$$\frac{(1-r)^{1+2\lambda}}{(1+r)^{n-1}}u(0)
\le u(r\zeta)\le
\frac{(1+r)^{1+2\lambda}}{(1-r)^{n-1}}u(0)$$
for $\g>-\frac{n}{2}$, and 

$$\frac{(1+r)^{1+2\lambda}}{(1-r)^{n-1}}u(0)
\le u(r\zeta)\le
\frac{(1-r)^{1+2\lambda}}{(1+r)^{n-1}}u(0)$$
for $\g<-\frac{n}{2}$.
\end{thm}
  
\vspace{3mm}\noindent{\bf Remark.} ~ Case $\lambda=0$ is the
classical Harnack Inequality in $B^n$.

%
%
          %
          %
%
\vspace{3mm}
\begin{cor} 
  Let $U$ be the potential function defined in $B^n$ by a positive Borel
  measure $\mu$ on $S^{n-1}$ as follows:
$$U(x)=\int_{S^{n-1}}\frac{1}{|x-\eta|^{n+2\lambda}}d\mu(\eta).$$

\vspace{2mm}\noindent For $\zeta\in S^{n-1}$, if $\g> -\frac{n}{2}$
(if $\g<-\frac{n}{2})$, the function
$$(1-r)^{n+2\lambda}U(r\zeta)$$
is decreasing (increasing) for $0\leq r<1$.
\end{cor}
%

          %
          %
%
\vspace{3mm} In Theorem 1.1, $\lambda=\frac{n}{2}-1$ corresponds to
the Laplace-Beltrami operator $\Delta_{n/2-1}$ and the Poincar\'e
metric.  It is known(\cite{brelot}) that given a positive invariant
harmonic function (solutions of $\Delta_{n/2-1}u=0$), there exists a
positive Borel measure $\mu$ on $S^{n-1}$, such that
$$u(x)=\int_{S^{n-1}}P_{n/2-1}(x,\zeta)d\mu(\zeta)$$
In this case the monotonicity property in Theorem
1.1 implies the following corollary.

\vspace{3mm}
\begin{cor} 
Let $u$ be a positive solution of $\Delta_{n/2-1}u=0$ in $B^n$.  Then
$$\left(\frac{1-r}{1+r}\right)^{n-1} u(r\zeta)\text{ ~decreasing in }r,$$
$$ \left(\frac{1+r}{1-r}\right)^{n-1} u(r\zeta)\text{ ~increasing in }r.$$

\end{cor}
%

          %
          %
%
\vspace{3mm}
\begin{cor} 
Let $u$ be a positive harmonic function with respect to the Laplace
operator ($\lambda=0$) defined in $B^n$ by a positive Borel measure
$\mu$ on $S^{n-1}$ with the Poisson Kernel $P_0$.  For $\zeta\in
S^{n-1},~0\leq r<1$, the function  
$$\dfrac{(1-r)^{n-1}}{1+r}u(r\zeta)$$
is decreasing and the function
$$\dfrac{(1+r)^{n-1}}{1-r}u(r\zeta)$$
is increasing.  In addition,
$$\lim_{r\to 1}(1-r)^{n-1}u(r\zeta)=2\mu(\{\zeta\}),$$
$$\lim_{r\to
  1}\dfrac{u(r\zeta)}{1-r} =
\int_{S^{n-1}}\dfrac{2}{|\zeta-\xi|^{n}}d\mu(\xi).$$
\end{cor}
Corollary 1.5 is the same as a result in \cite{pan}.

          %
          %
%
\vspace{3mm}
\begin{cor} 

  Let $B^n(R)$ be the open ball of radius $R$. Let $u(z)$ be an
  invariant harmonic function in $B^n(R)$ ({\sl a.k.a} $u(Rz)$ is
  invariant harmonic in $B^n$) defined by the Poisson kernel
  $P_{\lambda}(\frac{x}{R},\zeta)$.  For $\zeta\in S^{n-1}$, if $\lambda
  > -\frac{n}{2}$ (if $\lambda<-\frac{n}{2}$), the function
$$\frac{1}{R^{n-2-2\lambda}}\frac{(R-r)^{n-1}}{(R+r)^{1+2\lambda}}u(r\zeta)$$
is decreasing (increasing) and the function 
$$\frac{1}{R^{n-2-2\lambda}}\frac{(R+r)^{n-1}}{(R-r)^{1+2\lambda}}u(r\zeta)$$
is increasing (decreasing) in $r$ for $0\le r<R$.  The case
$\lambda=0$ gives the monotonicity of functions
$$\left(\frac{R-r}{R}\right)^{n-2}\frac{R-r}{R+r}u(r\zeta)
\qquad \text{and} \qquad
\left(\frac{R+r}{R}\right)^{n-2}\frac{R+r}{R-r}u(r\zeta),
$$
which implies that, $\forall x\in B^n(r), 0\le r<R$,
$$\left(\frac{R}{R+r}\right)^{n-2}\frac{R-r}{R+r}u(0)
\le u(x)
\le \left(\frac{R}{R-r}\right)^{n-2}\frac{R+r}{R-r}u(0)$$
--- the classical Harnack Inequality.  

\end{cor}
%

          %
          %

%
\vspace{3mm}
\begin{cor}  
  Let $u$ be a positive invariant harmonic function defined in $B^n$ by a
  positive Borel measure $\mu$ on $S^{n-1}$ with the Poisson kernel
  $P_\lambda$.  Let $~0\le r'\le r<1$.

\vspace{2mm}\noindent If $\g>-\frac{n}{2}$,
$$\dfrac{(1-r)^{n-1}}{(1+r)^{1+2\g}}\max_{|x|=r} u(x)
\le 
\dfrac{(1-r')^{n-1}}{(1+r')^{1+2\g}}\max_{|x|=r'} u(x)
$$
$$\dfrac{(1+r)^{n-1}}{(1-r)^{1+2\g}}\min_{|x|=r} u(x)
\ge 
\dfrac{(1+r')^{n-1}}{(1-r')^{1+2\g}}\min_{|x|=r'} u(x)
$$

\vspace{2mm}\noindent If $\g<-\frac{n}{2}$,
$$\dfrac{(1+r)^{n-1}}{(1-r)^{1+2\g}}\max_{|x|=r} u(x)
\le 
\dfrac{(1+r')^{n-1}}{(1-r')^{1+2\g}}\max_{|x|=r'} u(x)
$$
$$\dfrac{(1-r)^{n-1}}{(1+r)^{1+2\g}}\min_{|x|=r} u(x)
\ge 
\dfrac{(1-r')^{n-1}}{(1+r')^{1+2\g}}\min_{|x|=r'} u(x)
$$
\end{cor}
%

          %
          %

\vspace{9mm} Similar results are obtained in complex space
$\co^n$.  Let
$$P_\alpha(z,\zeta) =
\dfrac{(1-|z|^2)^{n+2\alpha}}{|1-z\cdot\ol\zeta|^{2n+2\alpha}},
~\alpha\in\mathbb{R}$$
be the Poisson-Szeg\"o kernel for the operator
$$\Delta_{\alpha,\beta} = 
4(1-|z|^2)\left\{ 
\sum_{i,j}(\delta_{i,j}-z_i\ol z_j)\frac{\p^2}{\p z_i\p \ol z_j} 
+\alpha\sum_{j}z_j\frac{\p}{\p z_j} 
+\beta\sum_{j}\ol z_j\frac{\p}{\p\ol z_j} - \alpha\beta\right\}$$
with $\alpha=\beta$, where $z\cdot\ol\zeta=\sum_{i=1}^n z_i\ol\zeta_i$.
Define
$$u(z)=\int_{S^{n-1}}P_\alpha(z,\zeta)d\mu(\zeta),
\quad \alpha\in\mathbb{R}.$$
\vspace{2mm}
\begin{thm} 
  Let $u$ be a positive invariant harmonic function defined in the unit
  ball $B^n\subset \mathbb C^n$ by a positive Borel measure $\mu$ on
  $S^{n-1}=\partial B^n$ with the Poisson-Szeg\"o kernel. Given $\zeta\in
  S^{n-1}$, if $\alpha> -n$ (if $\alpha<-n$), the function
$$\dfrac{(1-r)^{n}}{(1+r)^{n+2\alpha}}u(r\zeta)$$
is decreasing (increasing) for $0\leq r<1$, and the function
$$\dfrac{(1+r)^{n}}{(1-r)^{n+2\alpha}}u(r\zeta)$$
is increasing (decreasing) for $0\leq r<1$.  Also
$$\lim_{r\to 1}(1-r)^{n}u(r\zeta)=
\begin{cases}
2^{n+2\alpha}\mu(\{\zeta\}), & \alpha>-n\\
\infty, & \alpha<-n, ~ \mu(\{\zeta\}^c)>0\\
2^{n+2\alpha}\mu(\{\zeta\}), 
& \alpha<-n, ~ \mu(\{\zeta\}^c)=0
\end{cases}$$
and
$$\lim_{r\to 1}\dfrac{u(r\zeta)}{(1-r)^{n+2\alpha}} =
\int_{S^{n-1}}\dfrac{2^{n+2\alpha}}{|\zeta-\eta|^{2n+2\alpha}}d\mu(\eta).$$
\end{thm}

\vspace{2mm} The following theorem describes invariant harmonic
functions on the rays.

          %
          %
%
\vspace{2mm}
\begin{thm}  
  Let $u$ be a positive invariant harmonic function defined in the
  unit ball $B^n\subset \mathbb C^n$ by a positive Borel measure $\mu$
  on $S^{n-1}$ with the Poisson-Szeg\"o kernel.  Let $\zeta\in
  S^{n-1}$ and $~0\le r'\le r<1$.

\vspace{2mm}\noindent if $\alpha> -n$, 
\begin{eqnarray*}
\left(\!\frac{1\!-r}
{1\!-r'}\!\right)^{\!n+2\alpha}
\left(\!\!\frac{1\!+\!r'}{1\!+\!r}\right)^{\!n} 
\!u(r'\zeta)
\le u(r\zeta) \le
\left(\!\frac{1\!+r}
{1\!+r'}\!\right)^{\!n+2\alpha}
\left(\!\!\frac{1\!-\!r'}{1\!-\!r}\right)^{\!n} 
\!u(r'\zeta)
\end{eqnarray*}

\vspace{2mm}\noindent If $\alpha<-n$,

\begin{eqnarray*}
\left(\!\!\frac{1\!+r}{
1\!+r'}\right)^{\!\!-2n-2\alpha}
\left(\!\frac{1\!-\!r^2}{1\!-\!r'^2}\!\right)^{\!n+2\alpha} 
\!u(r'\zeta)
\le u(r\zeta) \le
\left(\!\frac{1\!-r}
{1\!-r'}\!\right)^{\!n+2\alpha}
\left(\!\!\frac{1\!+\!r'}{1\!+\!r}\right)^{\!n} 
\!u(r'\zeta)
\end{eqnarray*}

\vspace{2mm}\noindent For $r'=0$, the above becomes
$$\frac{(1-r)^{n+2\alpha}}{(1+r)^{n}}u(0)
\le u(r\zeta)\le
\frac{(1+r)^{n+2\alpha}}{(1-r)^{n}}u(0)$$
for $\alpha>-n$, and 
$$\frac{(1+r)^{n+2\alpha}}{(1-r)^{n}}u(0)
\le u(r\zeta)\le
\frac{(1-r)^{n+2\alpha}}{(1+r)^{n}}u(0)$$
for $\alpha <-n$.  
\end{thm}
\vspace{3mm}
%
%
\section{Proofs of Theorem 1.1 and its corollaries}
%
%

\vspace{3mm} We need the following two lemmas for the proof of Theorem
1.1.

          %
          %
\vspace{3mm}
\begin{lem} 
  Let $x\in\mathbb{R}^{n}$, $|x|=r$, $\zeta\in S^{n-1}$.  

\vspace{2mm}If $\lambda>-\frac{n}{2}$ then
{\small
\begin{equation}
-\frac{(n\!+\!2\lambda\!-\!(n\!-\!2\lambda\!\!-\!2)r)
(1\!\!-\!r^2)^{2\lambda}}{|x-\zeta|^{\!n+2\lambda}}
\!\!\leq\!
\frac{\p}{\p r}\frac{(1\!\!-\!r^2)^{1+2\lambda}}{|x-\zeta|^{n+2\lambda}}
\!\!\leq\!
\frac{(n\!+\!2\lambda\!+\!(n\!-\!2\lambda\!-\!2)r)
(1\!\!-\!r^2)^{\!2\lambda}}
{|x-\zeta|^{n+2\lambda}}
\label{lemma21top}
\end{equation}
}
If $\lambda<-\frac{n}{2}$, then
{\small
\begin{equation}
\frac{(n\!+\!2\lambda\!+\!(n\!-\!2\lambda\!\!-\!2)r)
(1\!\!-\!r^2)^{2\lambda}}{|x-\zeta|^{\!n+2\lambda}}
\!\!\leq\!
\frac{\p}{\p r}\frac{(1\!\!-\!r^2)^{1+2\lambda}}{|x-\zeta|^{n+2\lambda}}
\!\!\leq\!
-\frac{(n\!+\!2\lambda\!-\!(n\!-\!2\lambda\!-\!2)r)
(1\!\!-\!r^2)^{\!2\lambda}}
{|x-\zeta|^{n+2\lambda}}
\label{lemma21bottom}
\end{equation}
}
\end{lem}

          %
          %
%
\vspace{3mm}
\begin{proof} 
  ~Write $x=|x|\eta=r\eta, ~\eta\cdot\zeta=\sum_{i=1}^n \eta_i
  \zeta_i$.
$$\frac{\p}{\p r}|x-\zeta|^2=\frac{\p}{\p r}(|x|^2-2r\eta\cdot\zeta+1)
= 2(r-\eta\cdot\zeta),$$
then
\begin{equation*}
\begin{split}
  \frac{\p}{\p r}|x-\zeta|^{n+2\lambda} 
  & = \frac{\p}{\p r}(|x-\zeta|^2)^{\frac{n+2\lambda}{2}} \\
  & = \frac{n+2\lambda}{2}(|x-\zeta|^2)^{\frac{n+2\lambda}{2}-1}
  \frac{\p}{\p r}|x-\zeta|^2 \\
  &= (n+2\lambda)~|x-\zeta|^{n+2\lambda-2}~(r-\eta\cdot\zeta),
\end{split}
\end{equation*}
and
{\small
\begin{equation}
\begin{split}
  &\frac{\p}{\p r}\frac{(1-r^2)^{1+2\lambda}}{|x-\zeta|^{n+2\lambda}}\\
  =& \frac{(1+2\lambda)(1-r^2)^{2\lambda}(-2r)|x-\zeta|^{n+2\lambda} -
    (1-r^2)^{1+2\lambda}\frac{\p}{\p r}|x-\zeta|^{n+2\lambda}}
  {|x-\zeta|^{2(n+2\lambda)}}\\
  =&\frac{-2(1+2\lambda)(1\!\!-\!r^2)^{2\lambda}r|x\!-\!\zeta|^{n+2\lambda}
    - \!(1\!\!-\!r^2)^{1+2\lambda}(n+2\lambda)|x-\zeta|^{n+2\lambda-2}
    (r\!-\!\eta\cdot\zeta)}
  {|x-\zeta|^{2(n+2\lambda)}}\\
  =&\frac{-2(1+2\lambda)(1-r^2)^{2\lambda}r|x-\zeta|^{2}-(1-r^2)^{1+2\lambda}
    (n+2\lambda)(r-\eta\cdot\zeta)}{|x-\zeta|^{n+2\lambda+2}}.
\end{split}
\label{lemma21middle}
\end{equation}
}
To prove the right side inequality in (\ref{lemma21top}), it suffices to
show
$$
-2(1+2\lambda)r|x-\zeta|^{2} - (1-r^2)(n+2\lambda)(r-\eta\cdot\zeta)
\leq 
(n+2\lambda+(n-2\lambda-2)r)|x-\zeta|^2,
$$
which is equivalent to
$$-(n+2\lambda)(1-r^2)(r-\eta\cdot\zeta)\leq
(n+2\lambda)(1+r)|x-\zeta|^2.$$
For $\lambda>-\frac{n}{2}$, the above becomes
$$-(1-r^2)(r-\eta\cdot\zeta)\leq (1+r)|x-\zeta|^2,$$
or
$$-(1-r)(r-\eta\cdot\zeta)\leq r^2 -2\eta\cdot\zeta +1,$$
which, after a simple simplification, is equivalent to
$$\eta\cdot\zeta\leq 1$$
The inequality is true since $\zeta, \eta\in S^{n-1}$.  To prove the left
side inequality in (\ref{lemma21top}), it
suffices to show (using the result of (\ref{lemma21middle}))
$$-2(1+2\lambda)r|x-\zeta|^{2}-(1-r^2)(n+2\lambda)(r-\eta\cdot\zeta)
\geq -(n+2\lambda-(n-2\lambda-2)r)|x-\zeta|^2,$$
which is equivalent to
$$(n+2\lambda)(1-r^2)(r-\eta\cdot\zeta)\leq
(n+2\lambda)(1-r)|x-\zeta|^2.$$
For $\lambda>-\frac{n}{2}$, the inequality is equivalent to
$$(1-r^2)(r-\eta\cdot\zeta)\leq (1-r)|x-\zeta|^2,$$
which is, after a simplification, 
$$-\eta\cdot\zeta\leq 1,$$
true since $\zeta, \eta\in S^{n-1}$.  The proof of (\ref{lemma21bottom})
for $\lambda<-\frac{n}{2}$ is parallel.  This completes the proof of
Lemma 2.1.
\end{proof}
%

          %
          %
\vspace{2mm}
\begin{lem} 
  Let $u$ be a positive invariant harmonic function in $B^n$ defined by a
  positive Borel measure on $S^{n-1}$ with the Poisson kernel. 
\\

\noindent If $\lambda > -\frac{n}{2}$,
\begin{equation}
-\frac{(n+2\lambda-(n-2\lambda-2)r)}{1-r^2}u(x)
\leq
\frac{\partial u(x)}{\partial r}
\leq 
\frac{(n+2\lambda+(n-2\lambda-2)r)}{1-r^2}u(x).
\label{lemma22top}
\end{equation}
If $\lambda < -\frac{n}{2}$,
\begin{equation}
\frac{(n+2\lambda+(n-2\lambda-2)r)}{1-r^2}u(x)
\leq
\frac{\partial u(x)}{\partial r}
\leq -\frac{(n+2\lambda-(n-2\lambda-2)r)}{1-r^2}u(x).
\label{lemma22bottom}
\end{equation}

\end{lem}
%

          %
          %
%
\vspace{2mm}
\begin{proof} 
By the Poisson integral representation of $u$ in $B^n$,
$$u(x)=\int_{S^{n-1}}\frac{(1-|x|^2)^{1+2\g}}{|x-\zeta|^{n+2\g}}
d\mu(\zeta)$$
for a positive Borel measure $\mu$. By (\ref{lemma21top}) in Lemma 2.1 and
$\mu$ being a positive measure,
\begin{equation*}
\begin{split}
\int_{S^{n-1}}\frac{\partial}{\partial r}
\left(\frac{(1-|x|^2)^{1+2\g}}{|x-\zeta|^{n+2\g}}\right)d\mu(\zeta)
&\le\int_{S^{n-1}}\frac{(n+2\lambda+(n-2\lambda-2)r)(1-r^2)^{2\lambda}}
{|x-\zeta|^{n+2\lambda}}d\mu(\zeta)\\
&=\frac{(n+2\lambda+(n-2\lambda-2)r)}{1-r^2}
\int_{S^{n-1}}\frac{(1-|x|^2)^{1+2\g}}
{|x-\zeta|^{n+2\g}}d\mu(\zeta)\\
&=\frac{(n+2\lambda+(n-2\lambda-2)r)}{1-r^2}~u(x)
\end{split}
\end{equation*}
when $\lambda>-\frac{n}{2}$.  It follows that
$$\frac{\partial u(x)}{\partial r}
=\int_{S^{n-1}}\frac{\partial}{\partial r}
\left(\frac{(1-|x|^2)^{1+2\g}}{|x-\zeta|^{n+2\g}}\right)d\mu(\zeta)
\le\frac{(n+2\lambda+(n-2\lambda-2)r)}{1-r^2}~u(x).
$$
The left side inequality in (\ref{lemma22top}) can be proved in the same
manner.  For the equality case, consider
$u_y(x)=u(x,y)=\dfrac{(1-|x|^{2})^{1+2\lambda}}{|x-y|^{n+2\lambda}}$ which
is invariant harmonic in $\mathbb{R}^{n}\setminus\{y\}$ for $y\in
S^{n-1}$. A simple calculation shows that the equalities hold for $u_y(x)$
when $x=|x|y$ and $x=-|x|y$ respectively.  The proof of
(\ref{lemma22bottom}) is similar.  This completes the proof of Lemma 2.2.
\end{proof}
%

          %
          %

%
\vspace{3mm}Now we prove Theorem 1.1.

\vspace{2mm}\begin{proof} 
Consider $\varphi(r)=\dfrac{(1-r)^{n-1}}{(1+r)^{1+2\lambda}}$ and
$\psi(r)=\dfrac{(1+r)^{n-1}}{(1-r)^{1+2\lambda}}$ for $0\leq r <1$.
$$\frac{\varphi'}{\varphi}
=-\frac{(n+2\lambda+(n-2\lambda-2)r)}{1-r^2}~,$$
$$\frac{\psi'}{\psi}
=\frac{(n+2\lambda-(n-2\lambda-2)r)}{1-r^2}~.$$
Given $\omega\in
S^{n-1}$, consider
\begin{eqnarray*}
I(r,\omega)&=&\varphi(r)u(r\omega),\\
J(r,\omega)&=&\psi(r)u(r\omega).
\end{eqnarray*}
To show Theorem 1.1, it suffices to show that $I(r,\omega)$ is decreasing
(increasing) and $J(r,\omega)$ is increasing (decreasing) in $r$ for $0\le
r<1$ when $\lambda>-\frac{n}{2}$ (when $\lambda<-\frac{n}{2}$). By
(\ref{lemma22top}) in Lemma 2.2, for $\lambda>-\frac{n}{2}$,
\begin{equation*}
\begin{split}
\frac{d}{dr}(\log I(r,\omega))
&=\frac{\varphi'}{\varphi}+\frac{u'_r}{u}\\
&=-\frac{(n+2\lambda+(n-2\lambda-2)r)}{1-r^2}+\frac{u'_r}{u}\\
&\leq-\frac{(n+2\lambda+(n-2\lambda-2)r)}{1-r^2}
+\frac{(n+2\lambda+(n-2\lambda-2)r)}{1-r^2}\\
&=0.
\end{split}
\end{equation*}
Therefore $\log I(r,\omega)$ is decreasing in $r$, and so is $I(r,\omega)$.
Similarly,
\begin{equation*}
\begin{split}
\frac{d}{dr}(\log J(r,\omega))
&=\frac{\psi'}{\psi}+\frac{u'_r}{u}\\
&=\frac{(n+2\lambda-(n-2\lambda-2)r)}{1-r^2}+\frac{u'_r}{u}\\
&\geq\frac{(n+2\lambda-(n-2\lambda-2)r)}{1-r^2}
-\frac{(n+2\lambda-(n-2\lambda-2)r)}{1-r^2}\\
&=0.
\end{split}
\end{equation*}
Hence, $J(r,\omega)$ is increasing in $r$.  
%
%
For $\lambda>-\frac{n}{2}$ and $y\in S^{n-1}$,
\begin{eqnarray*}
\lim_{r\to 1}\frac{(1-r)^{n-1}}{(1+r)^{1+2\lambda}}P_\lambda(r\zeta,y)
&=&\lim_{r\to 1}\frac{(1-r)^{n-1}}{(1+r)^{1+2\lambda}}
\frac{(1-|r|^2)^{1+2\lambda}}{|r\zeta-y|^{n+2\lambda}}\\
&=&\lim_{r\to 1}\frac{(1-r)^{n+2\lambda}}{|r\zeta-y|^{n+2\lambda}}
\searrow \delta(\zeta,y)=
\begin{cases} 1,& \zeta=y\\0,& \zeta\not= y\end{cases}
\end{eqnarray*}
by applying the monotonicity properties in Theorem 1.1 to
$u=P_\lambda$.  By Lebesgue's dominated convergence theorem,
\begin{eqnarray*}
\lim_{r\to 1}(1-r)^{n-1}u(r\zeta)
&=&
\lim_{r\to 1}(1-r)^{n-1}\int_{S^{n-1}}P_\lambda(r\zeta,\xi)d\mu(\xi)\\
&=&
\lim_{r\to 1}(1+r)^{1+2\lambda}\int_{S^{n-1}}
\lim_{r\to 1}\frac{(1-r)^{n-1}}{(1+r)^{1+2\lambda}}
P_\lambda(r\zeta,\xi)d\mu(\xi)\\
&=& 2^{1+2\lambda}\mu(\{\zeta\}).
\end{eqnarray*}
Similarly, $\displaystyle\frac{(1+r)^{n-1}}{(1-r)^{1+2\lambda}}
P_\lambda(r\zeta,y)=\frac{(1+r)^{n+2\lambda}}
{|r\zeta-\xi|^{n+2\lambda}}$ is increasing in $r$ for
$\lambda>-\frac{n}{2}$.  By Lebesgue's monotone convergence theorem,
\begin{eqnarray*}
\lim_{r\to 1} \frac{u(r\zeta)}{(1-r)^{1+2\lambda}}
&=&
\lim_{r\to 1}\frac{1}{(1-r)^{1+2\lambda}}
\int_{S^{n-1}}P_\lambda(r\zeta,\xi)d\mu(\xi)\\
&=&
\lim_{r\to 1}\frac{1}{(1+r)^{n-1}}\int_{S^{n-1}}
\lim_{r\to 1}\frac{(1+r)^{n-1}}{(1-r)^{1+2\lambda}}
P_\lambda(r\zeta,\xi)d\mu(\xi)\\
&=&
\frac{1}{2^{n-1}}\int_{S^{n-1}}
\lim_{r\to 1}\frac{(1+r)^{n+2\lambda}}
{|r\zeta-\xi|^{n+2\lambda}}
d\mu(\xi)\\
&=&
\displaystyle
\int_{S^{n-1}}\frac{2^{1+2\lambda}}{|\zeta-\xi|^{n+2\lambda}}d\mu(\xi).
\end{eqnarray*}
For $\lambda<-\frac{n}{2}$, the monotonicity of $I(r,\omega)$ and
$J(r,\omega)$ is proved similarly to the case $\lambda>-\frac{n}{2}$ using
(\ref{lemma22bottom}) instead of (\ref{lemma22top}) in Lemma 2.2. 
$$
\lim_{r\to 1}\frac{(1-r)^{n-1}}{(1+r)^{1+2\lambda}}P_\lambda(r\zeta,y)
=\lim_{r\to 1}\frac{|r\zeta-y|^{-(n+2\lambda)}}{(1-r)^{-(n+2\lambda)}}
\nearrow
\begin{cases} 1,& \zeta=y\\\infty,& \zeta\not= y\end{cases}
$$
when $\lambda<-\frac{n}{2}$.  Therefore,
\begin{eqnarray*}
\lim_{r\to 1}(1-r)^{n-1}u(r\zeta)
&=&
\lim_{r\to 1}(1+r)^{1+2\lambda}\int_{S^{n-1}}
\lim_{r\to 1}\frac{(1-r)^{n-1}}{(1+r)^{1+2\lambda}}
P_\lambda(r\zeta,\xi)d\mu(\xi)\\
&=&
\begin{cases}2^
{1+2\lambda}\mu(\{\zeta\}),& if ~\mu(\{\zeta\}^c)=0;\\ \infty, & if
~\mu(\{\zeta\}^c)>0. 
\end{cases}
\end{eqnarray*}
Similarly, $\displaystyle
\frac{(1+r)^{n-1}}{(1-r)^{1+2\lambda}}P_\lambda(r\zeta,y)=
\frac{(1+r)^{n+2\lambda}} {|r\zeta-\xi|^{n+2\lambda}}$ is
decreasing in $r$ for $\lambda<-\frac{n}{2}$.  By Lebesgue's monotone
convergence theorem,
\begin{eqnarray*}
\lim_{r\to 1} \frac{u(r\zeta)}{(1-r)^{1+2\lambda}}
&=&
\lim_{r\to 1}\frac{1}{(1-r)^{1+2\lambda}}
\int_{S^{n-1}}P_\lambda(r\zeta,\xi)d\mu(\xi)\\
&=&
\lim_{r\to 1}\frac{1}{(1+r)^{n-1}}\int_{S_{n-1}}
\lim_{r\to 1}\frac{(1+r)^{n-1}}{(1-r)^{1+2\lambda}}
P_\lambda(r\zeta,\xi)d\mu(\xi)\\
&=&
\frac{1}{2^{n-1}}\int_{S^{n-1}}
\lim_{r\to 1}\frac{(1+r)^{n+2\lambda}}
{|r\zeta-\xi|^{n+2\lambda}}
d\mu(\xi)\\
&=&
\int_{S^{n-1}}\frac{2^{1+2\lambda}}{|\zeta-\xi|^{n+2\lambda}}d\mu(\xi)
\end{eqnarray*}
%
This completes the proof of Theorem 1.1.
\end{proof}
%

          %
          %
%
\vspace{3mm} The proof of Corollary 1.2 is straightforward and is omitted.
The proof of Corollary 1.3 follows.

\vspace{2mm}
\begin{proof} 
\begin{eqnarray*}
(1-r)^{n+2\lambda} U(r\zeta)
&=& \int_{S^{n-1}}\frac{(1-r)^{n+2\lambda}}
{|r\zeta-\eta|^{n+2\lambda}} d\mu(\eta) \\
&=& \frac{(1-r)^{n-1}}{(1+r)^{1+2\lambda}}
\int_{S^{n-1}}
\frac{(1-r^2)^{1+2\lambda}}{|r\zeta-\eta|^{n+2\lambda}}d\mu(\eta)\\
&=&\frac{(1-r)^{n-1}}{(1+r)^{1+2\lambda}} u(r\zeta)
\end{eqnarray*}
which is decreasing (increasing) in $r$ for $\lambda > - \frac{n}{2}$
($\lambda < -\frac{n}{2}$) by Theorem 1.1.  This completes the proof of
Corollary 1.3.
\end{proof}
%

          %
          %

%
\vspace{3mm} Corollaries 1.4 and 1.5 are special cases of Theorem 1.1.
Corollary 1.6 is a straightforward generalization from $B^n$ to $B^n(R)$.
The following is the proof of Corollary 1.7.

\vspace{3mm}\begin{proof}~ $0\le r'\le r<1$.  By the maximum
  principle, there is $\zeta\in S^{n-1}$ such that
  $u(r\zeta)=\max_{|x|=r}u(x)$.

\vspace{2mm}\noindent If $\g>-\frac{n}{2}$,  Theorem 1.1 implies 
\begin{eqnarray*}
\dfrac{(1-r)^{n-1}}{(1+r)^{1+2\g}}\max_{|x|=r} u(x)
&=&\dfrac{(1-r)^{n-1}}{(1+r)^{1+2\g}}u(r\zeta)\\
&\le& 
\dfrac{(1-r')^{n-1}}{(1+r')^{1+2\g}}u(r'\zeta)
\le
\dfrac{(1-r')^{n-1}}{(1+r')^{1+2\g}}\max_{|x|=r'} u(x)
\end{eqnarray*}
Similarly, there is $\xi\in S^{n-1}$ such that
$u(r\xi)=\min_{|x|=r}u(x)$.  When $\lambda>-\frac{n}{2}$, Theorem 1.1
yields
\begin{eqnarray*}
\dfrac{(1+r)^{n-1}}{(1-r)^{1+2\g}}\min_{|x|=r} u(x)
&=&\dfrac{(1+r)^{n-1}}{(1-r)^{1+2\g}}u(r\xi)\\
&\ge& 
\dfrac{(1+r')^{n-1}}{(1-r')^{1+2\g}}u(r'\xi)
\ge
\dfrac{(1+r')^{n-1}}{(1-r')^{1+2\g}}\min_{|x|=r'} u(x)
\end{eqnarray*}
The proof for $\lambda<-\frac{n}{2}$ is parallel.  This completes the
proof of Corollary 1.7.

\end{proof}

\vspace{3mm}
%
\section{Proof of Theorem 1.8}
%
%

\vspace{2mm} In the following, $B^n$ denotes the unit ball in
$\mathbb C^n$ and $S^{n-1}=\partial B^n$ the sphere.  We need
the following three lemmas for the proof of Theorem~1.8.

          %
          %
%
\vspace{2mm}
\begin{lem}
  If $a\in\mathbb{C}, |a|\leq 1$, then for $0\leq r\leq 1,$
\begin{equation}
1+r|a|^2\geq (1+r)\rp (a)
\label{lemma31top}
\end{equation}
\begin{equation}
1-r|a|^2\geq (-1+r)\rp (a)
\label{lemma31bottom}
\end{equation}
\end{lem}
%

          %
          %
%
\vspace{2mm}\begin{proof}
$|a|\le 1$, so $-1\le -|a|\le\rp(a)\le|a|\le 1$ and $\rp(a)^2\le
|a|^2$.  

\vspace{2mm} If $|a|^2\ge\rp(a)$, then $1+r|a|^2\ge \rp(a) + r\rp(a)$
so (\ref{lemma31top}) holds.  

\vspace{2mm} If $|a|^2<\rp(a)$, consider $f(r)=1+r|a|^2-(1+r)\rp(a),~
f'(r)=|a|^2-\rp(a)<0$.  So $f(r)$ decreases in $r\in[0,1]$.
$f(1)=1+|a|^2-2\rp(a)>1+\rp(a)^2-2\rp(a)=(1-Re(a))^2\ge 0$. So
$f(r)\ge 0$ and (\ref{lemma31top}) holds.  

\vspace{2mm} For the second inequality, $1+\rp(a)\ge
|a|^2+\rp(a)\ge r|a|^2+r\rp(a)$, so $1-r|a|^2\ge(-1+r)\rp(a)$ and
(\ref{lemma31bottom}) holds.
\end{proof}
%

          %
          %
%
\vspace{2mm}
\begin{lem} 
Let $z\in B^{n},~|z|=r,~\zeta\in S^{n-1}$.

\vspace{2mm}\noindent ~If $\alpha>-n$ then
{\small
\begin{equation}
-\frac{(2n\!+\!2\alpha+2\alpha r)(1\!-\!r^2)^{n+2\alpha\!-\!1}}
{|1\!-\!z\cdot\ol\zeta|^{n+2\alpha}}
\!\leq\!
\frac{\p}{\p r}\!
\frac{(1\!-\!r^2)^{n+2\alpha}}{|1\!-\!z\!\cdot\!\ol\zeta|^{2n+2\alpha}}
\!\leq\! \frac{(2n\!+\!2\alpha\!-\!2\alpha r)
(1\!-\!r^2)^{n+2\alpha\!-\!1}}{|1-z\cdot\ol\zeta|^{n+2\alpha}}
\label{lemma32top}
\end{equation}
}
If $\alpha<-n$, then
{\small
\begin{equation}
\frac{(2n\!+\!2\alpha\!-\!\!2\alpha r)(1\!-\!r^2)^{n+2\alpha\!-\!1}}
{|1-z\cdot\ol\zeta|^{n+2\alpha}}
\!\leq\!
\frac{\p}{\p
  r}\frac{(1-r^2)^{n+2\alpha}}{|1-z\!\cdot\!\ol\zeta|^{2n+2\alpha}}
\!\leq \!
- \frac{(2n\!+\!2\alpha\!+\!2\alpha r)(1\!-\!r^2)^{n+2\alpha\!-\!1}}
{|1-z\cdot\ol\zeta|^{n+2\alpha}}
\label{lemma32bottom}
\end{equation}
}
\end{lem}
%

          %
          %

\vspace{3mm}
\begin{proof}  
~ Let $z=|z|\eta = r\eta$.
$$\frac{\p}{\p r}|1-z\cdot\ol\zeta|^2
=\frac{\p}{\p r}(1-2r\rp (\eta\cdot\ol\zeta)
+ r^2|\eta\cdot\ol\zeta|^2)
=2(r|\eta\cdot\ol\zeta|^2-\rp(\eta\cdot\ol\zeta)),$$
and
\begin{equation*}
\begin{split}
\frac{\p}{\p r}|1-z\cdot\ol\zeta|^{2n+2\alpha} 
&=\frac{\p}{\p r}(|1-z\cdot\ol\zeta|^2)^{{n+\alpha}} \\
&=(n+\alpha)|(1-z\cdot\ol\zeta|^2)^{n+\alpha-1}
\frac{\p}{\p r}|1-z\cdot\ol\zeta|^2\\
&=2(n+\alpha)|1-z\cdot\ol\zeta|^{2n+2\alpha-2}
(r|\eta\cdot\ol\zeta|^2-\rp(\eta\cdot\ol\zeta)).
\end{split}
\end{equation*}
We have
{\small
\begin{equation}
\begin{split}
\frac{\p}{\p r}\frac{(1-r^2)^{n+2\alpha}}
{|1-z\cdot\ol\zeta|^{2n+2\alpha}}
&=\frac{(n+2\alpha)(1-r^2)^{n+2\alpha-1}
(-2r)|1-z\cdot\ol\zeta|^{2n+2\alpha}}{|1-z\cdot\ol\zeta|^{4n+4\alpha}}\\
&~~-~\frac{(1-r^2)^{n+2\alpha}
\frac{\p}{\p r}|1-z\cdot\ol\zeta|^{2n+2\alpha} }
{|1-z\cdot\ol\zeta|^{4n+4\alpha}}\\
&=\frac{-2(n+2\alpha)(1-r^2)^{n+2\alpha-1}r|1-z\cdot\ol\zeta|^{2n+2\alpha}}
{|1-z\cdot\ol\zeta|^{4n+4\alpha}}\\
&-
\frac{(1-r^2)^{n+2\alpha}2(n+\alpha)|1-z\cdot\ol\zeta|^{2n+2\alpha-2}
(r|\eta\cdot\ol\zeta|^2\!-\!\!\rp(\eta\cdot\ol\zeta))}
{|1-z\cdot\ol\zeta|^{4n+4\alpha}}\\
&=\frac{-2(n+2\alpha)(1-r^2)^{n+2\alpha-1}r|1-z\cdot\ol\zeta|^{2}}
{|1-z\cdot\ol\zeta|^{2n+2\alpha+2}}\\
&~-\frac{(1-r^2)^{n+2\alpha}2(n+\alpha)
(r|\eta\cdot\ol\zeta|^2-\rp(\eta\cdot\ol\zeta))}
{|1-z\cdot\ol\zeta|^{2n+2\alpha+2}}
\end{split}
\end{equation}
}
To prove the right side inequality of (\ref{lemma32top}), it suffices
to prove
$$-2(n+2\alpha)r|1-z\cdot\ol\zeta|^{2}-(1-r^2)
(2n+2\alpha)(r|\eta\cdot\ol\zeta|^2-\rp(\eta\cdot\ol\zeta))
\leq (2n+2\alpha-2\alpha r)|1-z\cdot\ol\zeta|^2$$
which is equivalent to
$$-(1-r^2)2(n+\alpha)(r|\eta\cdot\ol\zeta|^2-\rp(\eta\cdot\ol\zeta))
\leq 2(n+\alpha)(1+r)|1-z\cdot\ol\zeta|^2. $$
For $\alpha>-n$, the above inequality is equivalent to
$$-(1-r)((r|\eta\cdot\ol\zeta|^2-\rp(\eta\cdot\ol\zeta))\leq
|1-z\cdot\ol\zeta|^2$$
which is, after a simple simplification, 
$$(1+r)\rp(\eta\cdot\ol\zeta)\leq 1+r|\eta\cdot\ol\zeta|^2.$$
The inequality is true by (\ref{lemma31top}) in Lemma 3.1.  To prove the
left side inequality of (\ref{lemma32top}), it suffices to show
$$-2(n+2\alpha)r|1-z\cdot\ol\zeta|^{2}-(1-r^2)
(2n+2\alpha)(r|\eta\cdot\ol\zeta|^2-\rp(\eta\cdot\ol\zeta))
\geq -(2n+2\alpha+2\alpha r)|1-z\cdot\ol\zeta|^2$$
which is equivalent to
$$-(1-r^2)2(n+\alpha)(r|\eta\cdot\ol\zeta|^2-\rp(\eta\cdot\ol\zeta))
\geq -2(n+\alpha)(1-r)|1-z\cdot\ol\zeta|^2.$$
For $\alpha>-n$, the above inequality becomes
$$(1+r)((r|\eta\cdot\ol\zeta|^2-\rp(\eta\cdot\ol\zeta))
\leq |1-z\cdot\ol\zeta|^2$$
which is, after a simple simplification, 
$$(-1+r)\rp(\eta\cdot\ol\zeta)\leq 1-r|\eta\cdot\ol\zeta|^2.$$
The inequality is true by (\ref{lemma31bottom}) in Lemma 3.1.  The
proof of (\ref{lemma32bottom}) is parallel to that of
(\ref{lemma32top}), using the same inequalities in Lemma 3.1.
This completes the proof of Lemma 3.2.
\end{proof}
%

          %
          %
%
\vspace{2mm}\begin{lem}  
  Let $u(z)$ be a positive invariant harmonic function in $B^n$ defined by
  a positive Borel measure on $S^{n-1}$ with the Poisson-Szeg\"o kernel,
  $|z|=r$.  

\noindent If $\alpha>-n$,
\begin{equation}
-\frac{(2n+2\alpha+2\alpha r)}{1-r^2}u(z)
\leq\frac{\partial u(z)}{\partial r}
\leq \frac{(2n+2\alpha-2\alpha r)}{1-r^2}u(z).
\label{lemma33top}
\end{equation}
If $\alpha<-n$,
\begin{equation}
\frac{(2n+2\alpha-2\alpha r)}{1-r^2}u(z)
\leq\frac{\partial u(z)}{\partial r}
\leq -\frac{(2n+2\alpha+2\alpha r)}{1-r^2}u(z).
\label{lemma33bottom}
\end{equation}
\end{lem}
%

          %
          %
\vspace{2mm}
\begin{proof} 
By the Poisson-Szeg\"o integral representation of $u$ in $B^n$,
$$u(z)=\int_{S^{n-1}}
\frac{(1-|z|^2)^{n+2\alpha}}{|1-z\cdot\ol\zeta|^{2n+2\alpha}}d\mu(\zeta)$$
for a positive Borel measure $\mu$ on $S^{n-1}$.  By (\ref{lemma32top}) in
Lemma 3.2 and $\mu$ being a positive measure,
\begin{equation*}
\begin{split}
\int_{S^{n-1}}\frac{\partial}{\partial r}
\left(\frac{(1-|z|^2)^{n+2\alpha}}{|1-z\cdot\ol\zeta|^{2n+2\alpha}}\right)
d\mu(\zeta)
&\leq\int_{S^{n-1}}
\frac{(n+2\alpha-2\alpha r)(1-r^2)^{n+2\alpha-1}}
{|1-z\cdot\ol\zeta|^{n+2\alpha}}d\mu(\zeta)\\
&=\frac{(n+2\alpha-2\alpha r)}{1-r^2}
\int_{S^{n-1}}\frac{(1-|z|^2)^{n+2\alpha}}
{|1-z\cdot\ol\zeta|^{2n+2\alpha}}d\mu(\zeta)\\
&=\frac{(n+2\alpha-2\alpha r)}{1-r^2}u(z)
\end{split}
\end{equation*}
when $\alpha>-n$.  It follows that
$$\frac{\partial u(z)}{\partial r}
=\int_{S^{n-1}}\frac{\partial}{\partial r}
\left(\frac{(1-|z|^2)^{n+2\alpha}}{|1-z\cdot\ol\zeta|^{2n+2\alpha}}\right)
d\mu(\zeta)
\le \frac{(n+2\alpha-2\alpha r)}{1-r^2}u(z).$$
The left side inequality in (\ref{lemma33top}) is proved similarly.
For the equality case, consider
$u_w(z)=P_\alpha(z,w)=\dfrac{(1-|z|^{2})^{n+2\alpha}}{|z-w|^{2n+2\alpha}}$.
It is known that $u_w(z)$ is invariant harmonic in $\mathbb
C^{n}\setminus\{w\}$ for $w\in S^{n-1}$. A simple calculation shows
that the equalities in (\ref{lemma33top}) hold for $u_w(z)$ when
$z=|z|w$ and $z=-|z|w$ respectively.  The proof of
(\ref{lemma33bottom}) is parallel to that of (\ref{lemma33top}), using
(\ref{lemma32bottom}) instead of (\ref{lemma32top}) in Lemma 3.2.
This completes the proof of Lemma 3.3.
\end{proof}
%
          %
          %

\vspace{2mm} Now we prove Theorem 1.8.

\vspace{2mm}\begin{proof} 
Consider $\varphi(r)=\dfrac{(1-r)^{n}}{(1+r)^{n+2\alpha}},~
\psi(r)=\dfrac{(1+r)^{n}}{(1-r)^{n+2\alpha}}$ for $0\leq r <1$.
$$\frac{\varphi'}{\varphi}=
-\frac{2n+2\alpha-2\alpha r}{1-r^2},$$
$$\frac{\psi'}{\psi}=\frac{2n+2\alpha+2\alpha r}{1-r^2}.$$
Given $\omega\in S^{n-1}$, consider
$$I(r,\omega)=\varphi(r)u(r\omega),$$
$$J(r,\omega)=\psi(r)u(r\omega).$$
To show Theorem 1.8, it suffices to show that $I(r,\omega)$ is decreasing
(increasing) and $J(r,\omega)$ is increasing (decreasing) in $r$ when
$\alpha > -n$ (when $\alpha < -n$).  By (\ref{lemma33top}) in Lemma 3.3,
when $\alpha > -n$,
\begin{equation*}
\begin{split}
\frac{\partial}{\partial r}\log I(r,\omega)
&=\frac{\varphi'}{\varphi}+\frac{u'_r}{u}\\
&=-\frac{2n+2\alpha-2\alpha r}{1-r^2}+\frac{u'_r}{u}\\
&\leq  -\frac{2n+2\alpha-2\alpha r}{1-r^2} 
+\frac{2n+2\alpha-2\alpha r}{1-r^2}\\
&=0.
\end{split}
\end{equation*}
Therefore $\log I(r,\omega)$ is decreasing in $r$, and so is
$I(r,\omega)$.  Similarly,
\begin{equation*}
\begin{split}
\frac{\partial}{\partial r}\log J(r,\omega)
&=\frac{\psi'}{\psi}+\frac{u'_r}{u}\\
&=\frac{2n+2\alpha+2\alpha r}{1-r^2}+\frac{u'_r}{u}\\
&\geq \frac{2n+2\alpha+2\alpha r}{1-r^2}
-\frac{2n+2\alpha+2\alpha r}{1-r^2}\\
&=0.
\end{split}
\end{equation*}
Hence, $J(r,\omega)$ is increasing in $r$. 
%
%
For $\alpha>-n$ and $\zeta,w\in S^{n-1}$,
\begin{eqnarray*}
\lim_{r\to 1}\frac{(1-r)^{n}}{(1+r)^{n+2\alpha}}P_\alpha(r\zeta,w)
&=&\lim_{r\to 1}\frac{(1-r)^{n}}{(1+r)^{n+2\alpha}}
\frac{(1-|r|^2)^{n+2\alpha}}{|1-r\zeta\cdot\ol w|^{2n+2\alpha}}\\
&=&\lim_{r\to 1}\frac{(1-r)^{2n+2\alpha}}{|1- r\zeta\cdot\ol w|^{2n+2\alpha}}
\searrow \delta(\zeta,w)
\end{eqnarray*}
by applying the monotonicity results in Theorem 1.8 to $u=P_\alpha$.
By Lebesgue's dominated convergence theorem,
\begin{eqnarray*}
\lim_{r\to 1}(1-r)^{n}u(r\zeta)
&=&
\lim_{r\to 1}(1-r)^{n}\int_{S^{n-1}}P_\alpha(r\zeta,\xi)d\mu(\xi)\\
&=&
\lim_{r\to 1}(1+r)^{n+2\alpha}\int_{S^{n-1}}
\lim_{r\to 1}\frac{(1-r)^{n}}{(1+r)^{n+2\alpha}}
P_\lambda(r\zeta,\xi)d\mu(\xi)\\
&=& 2^{n+2\alpha}\mu(\{\zeta\}).
\end{eqnarray*}
Similarly, $\displaystyle
\frac{(1+r)^{n}}{(1-r)^{n+2\alpha}}P_\alpha(r\zeta,w)=\frac{(1+r)^{2n+2\alpha}}
{|1-r\zeta\cdot\ol w|^{2n+2\alpha}}$ is increasing
in $r$ for $\alpha>-n$.  By Lebesgue's monotone convergence theorem,
\begin{eqnarray*}
\lim_{r\to 1} \frac{u(r\zeta)}{(1-r)^{n+2\alpha}}
&=&
\lim_{r\to 1}\frac{1}{(1-r)^{n+2\alpha}}
\int_{S^{n-1}}P_\alpha(r\zeta,\xi)d\mu(\xi)\\
&=&
\lim_{r\to 1}\frac{1}{(1+r)^{n}}\int_{S^{n-1}}
\lim_{r\to 1}\frac{(1+r)^{n}}{(1-r)^{n+2\alpha}}
P_\alpha(r\zeta,\xi)d\mu(\xi)\\
&=&
\frac{1}{2^{n}}\int_{S^{n-1}}
\lim_{r\to 1}\frac{(1+r)^{2n+2\alpha}}
{|1-r\zeta\cdot\ol\xi|^{2n+2\alpha}}
d\mu(\xi)\\
&=&
\displaystyle
\int_{S^{n-1}}\frac{2^{n+2\alpha}}{|\zeta-\xi|^{2n+2\alpha}}d\mu(\xi).
\end{eqnarray*}
For $\alpha<-n$, the monotonicity of $I(r,\omega)$ and $J(r,\omega)$ is
proved similarly by applying (\ref{lemma33bottom}) in Lemma 3.3.  
$$
\lim_{r\to 1}\frac{(1-r)^{n}}{(1+r)^{n+2\alpha}}P_\alpha(r\zeta,w)
=\lim_{r\to 1}\frac{|1-r\zeta\cdot\ol w|^{-(2n+2\alpha)}}{(1-r)^{-(2n+2\alpha)}}
\nearrow
\begin{cases} 1,& \zeta=w \\ \infty,& \zeta\not= w \end{cases}
$$
when $\alpha<-n$.  Therefore,
\begin{eqnarray*}
\lim_{r\to 1}(1-r)^{n}u(r\zeta)
&=&
\lim_{r\to 1}(1+r)^{n+2\alpha}\int_{S^{n-1}}
\lim_{r\to 1}\frac{(1-r)^{n}}{(1+r)^{n+2\alpha}}
P_\alpha(r\zeta,\xi)d\mu(\xi)\\
&=&
\begin{cases}2^
{n+2\alpha}\mu(\{\zeta\}),& if ~\mu(\{\zeta\}^c)=0;\\ \infty, &  if
~\mu(\{\zeta\}^c)>0. 
\end{cases}
\end{eqnarray*}
Similarly, $\displaystyle
\frac{(1+r)^{n}}{(1-r)^{n+2\alpha}}P_\alpha(r\zeta,w)=\frac{(1+r)^{2n+2\alpha}}
{|1-r\zeta\cdot\ol w|^{2n+2\alpha}}$ is decreasing
in $r$ for $\alpha<-n$.  By Lebesgue's monotone convergence theorem,
\begin{eqnarray*}
\lim_{r\to 1} \frac{u(r\zeta)}{(1-r)^{n+2\alpha}}
&=&
\lim_{r\to 1}\frac{1}{(1-r)^{n+2\alpha}}
\int_{S^{n-1}}P_\alpha(r\zeta,\xi)d\mu(\xi)\\
&=&
\lim_{r\to 1}\frac{1}{(1+r)^{n}}\int_{S_{n-1}}
\lim_{r\to 1}\frac{(1+r)^{n}}{(1-r)^{n+2\alpha}}
P_\alpha(r\zeta,\xi)d\mu(\xi)\\
&=&
\frac{1}{2^{n}}\int_{S^{n-1}}
\lim_{r\to 1}\frac{(1+r)^{2n+2\alpha}}
{|1-r\zeta\cdot\ol\xi|^{2n+2\alpha}}
d\mu(\xi)\\
&=&
\int_{S^{n-1}}\frac{2^{n+2\alpha}}{|\zeta-\xi|^{2n+2\alpha}}d\mu(\xi)
\end{eqnarray*}
%
This completes the proof of Theorem 1.8.
\end{proof}
%

       %
       %
%
\vspace{3mm}\noindent{\bf Remark}.  Notice that the monotonicity of the
auxiliary functions $\varphi$ and $\psi$ in the proofs of Theorem 1.1 and
Theorem 1.8 may vary depending on the values of the parameter $\lambda$
(or $\alpha$) and the dimension $n$.  When $\lambda>-\frac{n}{2}$ ( or
$\alpha>-n$), we have $\varphi'<0$ and $\psi'>0$, i.e. $\varphi$ increases
and $\psi$ decreases in $r$ for $0<r<1$.  For $\lambda<-\frac{n}{2}$ (or
$\alpha<-n$), the monotonicity does not necessarily hold.  For example, in
the real case in Theorem 1.1, for $\lambda<-\frac{n}{2}$,
$$\varphi(r)=\dfrac{(1-r)^{n-1}}{(1+r)^{1+2\lambda}}~,\qquad
\varphi'(r)
\begin{cases}
>0, &  
r\in \left(0,\displaystyle\frac{-2\lambda-n}{-2\lambda+(n-2)}\right) \\
<0, &  
r\in\left(\displaystyle\frac{-2\lambda-n}{-2\lambda+(n-2)},1\right)
\end{cases}
$$
i.e. the monotonicity may change for certain combinations of $n$ and
$\lambda$.  However, the monotonicity of $\varphi(r)u(r\zeta)$ and
$\psi(r)u(r\zeta)$ holds.  

\vspace{5mm}
%
\section{Proofs of Theorem 1.2 and Theorem 1.9}
%
%

\vspace{3mm}\noindent The proofs for the two theorems are based on the
following lemma.
%

          %
          %
%

\vspace{3mm}
\begin{lem}  
  Let $f(r)$ be a positive function on $r\in[0,1)$.  If for
  $a,b\in\mathbb R$,
\begin{equation}
-\frac{a+br}{1-r^2}f(r) \le f'(r) \le\frac{a-br}{1-r^2}f(r), 
\label{lemma41a}
\end{equation}
then for $0\le r'\le r<1$, 
\begin{equation}
\left(\!\!\frac{1\!+r}{
1\!+r'}\right)^{-a}
\left(\!\frac{1\!-\!r^2}{1\!-\!r'^2}\!\right)^{\!\frac{b+a}{2}} f(r')
\le f(r)  \le
\left(\!\frac{1\!+r}
{1\!+r'}\!\right)^{a}
\left(\!\!\frac{1\!-\!r^2}{1\!-\!r'^2}\right)^{\!\frac{b-a}{2}} f(r').
\label{lemma41b}
\end{equation}
\end{lem}
%

          %
          %

\vspace{3mm}
\begin{proof} 
\begin{equation*}
\int\frac{a-br}{1-r^2}dr 
= a\ln(1+r)+\frac{1}{2}(b-a)\ln(1-r^2)+C
\end{equation*}
Thus for $0\le r'\le r''<1$, by (\ref{lemma41a}),
\begin{eqnarray*}
\ln f(r'')-\ln f(r')=\int_{r'}^{r''}\frac{f'(r)}{f(r)}dr \le
\int_{r'}^{r''}\frac{a-br}{1-r^2}dr 
\le \ln \left(\!\frac{1\!+r''}
{1\!+r'}\!\right)^{a}
\left(\!\!\frac{1\!-\!r''^2}{1\!-\!r'^2}\right)^{\!\frac{b-a}{2}} 
\end{eqnarray*}
i.e. the right side inequality in (\ref{lemma41b}) holds.  Similarly,
by the left side of (\ref{lemma41a}),
\begin{eqnarray*}
\ln f(r'')-\ln f(r') 
\ge
\int_{r'}^{r''}-\frac{a+br}{1-r^2}dr 
\ge \ln \left(\!\frac{1\!+r''}
{1\!+r'}\!\right)^{-a}
\left(\!\!\frac{1\!-\!r''^2}{1\!-\!r'^2}\right)^{\!\frac{b+a}{2}}. 
\end{eqnarray*}
i.e. the
left side inequality in (\ref{lemma41b}) holds.
\end{proof}
%

          %
          %
%

\vspace{3mm}\noindent Now we prove Theorem 1.2.

\vspace{3mm}\begin{proof}
  If $\lambda>-\frac{n}{2}$, $u(r\zeta)$ satisfies (\ref{lemma22top})
  in Lemma 2.2.  Therefore (\ref{lemma41a}) holds with
  $f(r)=u(r\zeta),~a=n+2\lambda, ~b= -n+2\lambda+2$. Let $0\le r'\le
  r<1$.  (\ref{lemma41b}) in Lemma 4.1 implies
\begin{equation*}
\left(\!\!\frac{1\!+r}{
1\!+r'}\right)^{-n-2\lambda}
\left(\!\frac{1\!-\!r^2}{1\!-\!r'^2}\!\right)^{\! 2\lambda+1} 
\!\!\!\!\!u(r'\zeta)
\le u(r\zeta) \le
\left(\!\frac{1\!+r}
{1\!+r'}\!\right)^{n+2\lambda}
\left(\!\!\frac{1\!-\!r^2}{1\!-\!r'^2}\right)^{\!-n+1} 
\!\!\!u(r'\zeta).
\end{equation*}
If $\lambda<-\frac{n}{2}$, $u(r\zeta)$ satisfies (\ref{lemma22bottom})
in Lemma 2.2, thus (\ref{lemma41a}) holds with
$f(r)=u(r\zeta),~a=-n-2\lambda, ~b=-n+2\lambda+2$. Applying
(\ref{lemma41b}),
\begin{equation*}
\left(\!\!\frac{1\!+r}{
1\!+r'}\right)^{n+2\lambda}
\left(\!\frac{1\!-\!r^2}{1\!-\!r'^2}\!\right)^{\!-n+1} 
\!\!\!\!\!u(r'\zeta)
\le u(r\zeta) \le
\left(\!\frac{1\!+r}
{1\!+r'}\!\right)^{-n-2\lambda}
\left(\!\!\frac{1\!-\!r^2}{1\!-\!r'^2}\right)^{\!2\lambda+1} 
\!\!\!u(r'\zeta).
\end{equation*}
This completes the proof of Theorem 1.2.
\end{proof}
%
%

          %
          %
%

\vspace{3mm}\noindent The proof of Theorem 1.9 is similar to that of
Theorem 1.2.

\vspace{3mm}\begin{proof}
  If $\alpha>-n$, $u(r\zeta)$ satisfies (\ref{lemma33top}) in Lemma
  3.3.  Therefore (\ref{lemma41a}) holds with
  $f(r)=u(r\zeta),~a=2n+2\alpha, ~b=2\alpha$. Let $0\le r'\le r<1$.
  (\ref{lemma41b}) in Lemma 4.1 implies
\begin{equation*}
\left(\!\!\frac{1\!+r}{
1\!+r'}\right)^{-2n-2\alpha}
\left(\!\frac{1\!-\!r^2}{1\!-\!r'^2}\!\right)^{\! n+2\alpha} 
\!\!\!\!\!u(r'\zeta)
\le u(r\zeta) \le
\left(\!\frac{1\!+r}
{1\!+r'}\!\right)^{2n+2\alpha}
\left(\!\!\frac{1\!-\!r^2}{1\!-\!r'^2}\right)^{\!-n} 
\!\!\!u(r'\zeta).
\end{equation*}
If $\alpha<-n$, $u(r\zeta)$ satisfies (\ref{lemma33bottom}) in Lemma
3.3, thus (\ref{lemma41a}) holds with $f(r)=u(r\zeta),~a=-2n-2\alpha,
~b=2\alpha$. From (\ref{lemma41b}),
\begin{equation*}
\left(\!\!\frac{1\!+r}{
1\!+r'}\right)^{2n+2\alpha}
\left(\!\frac{1\!-\!r^2}{1\!-\!r'^2}\!\right)^{\!-n} 
\!\!\!\!\!u(r'\zeta)
\le u(r\zeta) \le
\left(\!\frac{1\!+r}
{1\!+r'}\!\right)^{-2n-2\alpha}
\left(\!\!\frac{1\!-\!r^2}{1\!-\!r'^2}\right)^{\!n+2\alpha} 
\!\!\!u(r'\zeta).
\end{equation*}
This completes the proof of Theorem 1.9.
\end{proof}

\vspace{3mm}\noindent Most results in this paper are on the
function values at two points in $B^n$ on the same ray.  Similar
results can be obtained for any two points in $B^n$ (\cite{pan2}).

\vspace{9mm}

\end{document}